\documentclass[11pt]{article}
\usepackage{amsmath, amssymb, eucal, latexsym, multicol, amsthm, bm}
\usepackage{graphicx,psfrag,color}
\usepackage{comment}
\usepackage{soul}

\usepackage[figtopcap,FIGBOTCAP]{subfigure}

\usepackage{enumerate}

\newcommand{\rhn}[1]{{\color{black}{#1}}}

\newcommand{\cV}{\mathcal{V}}
\newcommand{\cW}{\mathcal{W}}

\DeclareMathOperator{\clos}{clos}

\title{{\bf A saturation property for the spectral-Galerkin approximation of a Dirichlet problem in a square}}
\author{C. Canuto, R.H. Nochetto, R. Stevenson, and M. Verani}

\date{\today} 

\newtheorem{theorem}{Theorem}[section]
\newtheorem{lemma}[theorem]{Lemma}
\newtheorem{proposition}[theorem]{Proposition}

\theoremstyle{definition}

\theoremstyle{remark}
\newtheorem{remark}[theorem]{Remark}

\numberwithin{equation}{section}

\usepackage{color}

\begin{document}

\maketitle

\begin{abstract}
Both practice and analysis of adaptive $p$-FEMs and $hp$-FEMs raise
the question what increment in the current polynomial degree $p$
guarantees a $p$-independent reduction of the Galerkin error. We
answer this question for the $p$-FEM in the simplified context of homogeneous
Dirichlet problems for the Poisson equation in the
two dimensional unit square with polynomial data of
degree $p$. We show that an increment proportional to $p$ yields a
$p$-robust error reduction and provide computational evidence that a
constant increment does not.  
\end{abstract}

\section{Motivation and statement of the result}\label{sec:intro}

High order finite element methods (FEMs) exhibit exponential
convergence for elliptic problems with piecewise analytic data, and
thus have become the methods of choice in computational science and
engineering for such problems. The seminal work of Babu\v ska and
collaborators \cite{Babuska-1991,GuoBabuska-1,GuoBabuska-2}
has established the mathematical foundations for
the a priori design of meshes and distribution of polynomial degrees,
and proved exponential convergence for corner and edge singularities.
In contrast, adaptive $hp$-FEMs hinge on a posteriori error
estimators, which help determine whether it is more convenient to
locally refine the mesh or increase the polynomial degree to improve
the resolution. Although exponential convergence is observe
experimentally, it \rhn{has never been} proved rigorously with the exception of
\cite{CNSV16}.

Our adaptive $hp$-FEM of \cite{CNSV16} hinges on a coarsening module due to
Binev \cite{Binev-2015}, which in turn guarantees instance optimality and
thus exponential convergence. As any other adaptive $hp$-FEM, ours
also has a module to reduce the PDE error by a fixed fraction for
piecewise polynomial data thereby avoiding data oscillation. Such
module in \cite{CNSV16} relies on the a posteriori error estimator of
Melenk and \rhn{Wohlmuth} \cite{MW01} for dimension $d=2$ and does not possess
optimal complexity. In \cite{CNSV17} we turn to the equilibrated flux
residual estimator of Braess, Pillwein and Sch{\"o}berl \cite{BPS-2009}, and
Ern and Vorahl\'\i k \cite{ErnVohralik-2015,ErnVohralik-2016},
and show that the issue of optimal complexity
reduces to studying three model problems with polynomial data in the
reference triangle for $d=2$. We present \rhn{overwhelming} computational
evidence in \cite{CNSV17} supporting the fact that to reduce the Galerkin
error by a fixed factor, the polynomial degree $p$ must be increased by
an amount proportional to $p$.

In this paper we take over this question again in a further simplified
setting and give a rigorous answer. We consider the Poisson equation
\begin{equation}\label{poisson}
  -\Delta u = f \quad\textrm{in }\Omega,
  \qquad
  u = 0 \quad\textrm{on }\partial\Omega,
\end{equation}
over the unit square $\Omega = (-1,1)^2$ of $\mathbb{R}^2$ with
polynomial $f$ of degree $p$.

Let us first introduce some notation. Let $I = (-1,1)$ be the reference
element so that $\Omega = I^2$. For $p\geq 0$, let ${\mathbb
  P}_p(\Omega)$ denote the space of polynomials of total degree $\leq
p$ restricted to $\Omega$,
and let $\cV_p := {\mathbb P}_p(\Omega) \cap H^1_0(\Omega)$.
Since the latter space \rhn{reduces to} $\{0\}$ for $p<4$, we will consider it only for $p\geq 4$.
We equip $H^1_0(\Omega)$ with the energy inner product $(\cdot,\cdot)_{H^1_0(\Omega)}:=(\nabla\cdot,\nabla\cdot)_{L^2(\Omega)^2}$ and resulting norm $\|\cdot\|_{H^1_0(\Omega)}$. Analogous definitions apply to $H^1_0(I)$.

Given any $f \in {\mathbb P}_p(\Omega)$, let $u=u(f) \in H^1_0(\Omega)$ be the variational solution of \eqref{poisson}, i.e.,
\begin{equation}\label{eq:def-u}
(u,v)_{H^1_0(\Omega)} = (f,v)_{L^2(\Omega)} \qquad \forall v \in H^1_0(\Omega).
\end{equation}
For any $q \geq 4$, let $u_q =u_q(f) \in \cV_q$ be the corresponding
Galerkin projection of $u$ \rhn{onto $\cV_q$}, i.e., 
\begin{equation}\label{eq:def-uq}
(u_q,v_q)_{H^1_0(\Omega)} = (f,v_q)_{L^2(\Omega)} \qquad \forall v \in \cV_q.
\end{equation}

We are interested in finding sufficient conditions on $q=q(p)>p$ so that
\begin{equation}\label{error-reduction}
\|u-u_q\|_{H^1_0(\Omega)} \le \alpha \|u-u_p\|_{H^1_0(\Omega)}
\end{equation}
for $0<\alpha<1$ independent of $p$. If \eqref{error-reduction} holds,
we say that the error reduction is {\it $p$-robust.} In view of
Phytagoras equality
\begin{equation}\label{phytagoras}
  \|u-u_p\|_{H^1_0(\Omega)}^2 = \|u-u_q\|_{H^1_0(\Omega)}^2
  + \|u_q-u_p\|_{H^1_0(\Omega)}^2,
\end{equation}
which is a consequence of Galerkin orthogonality $u_q-u_p \perp u-u_q$
in $H^1_0(\Omega)$, we see that \eqref{error-reduction} is equivalent
to the following {\it saturation property}
\begin{equation}\label{saturation-1}
\sqrt{1-\alpha^2}\, \|u-u_p\|_{H^1_0(\Omega)} \le \|u_q-u_p\|_{H^1_0(\Omega)}.
\end{equation}
%
%
%
We next
observe that \eqref{saturation-1} is equivalent to the simpler
saturation property
\begin{equation}\label{saturation-2}
\sqrt{1-\alpha^2}\,  \|u\|_{H^1_0(\Omega)} \le \|u_q\|_{H^1_0(\Omega)}.
\end{equation}
To see this just define $v := u - u_p \in H^1_0(\Omega)$, which is the
solution of \eqref{poisson} with polynomial forcing
$f-\Delta u_p \in \rhn{\mathbb{P}_p(\Omega)}$ and Galerkin solutions $v_p=0$ and $v_q = u_q-u_p$.
We aim at establishing the following rigorous result.

\begin{theorem}[saturation property]\label{theo:main}
There exists a constant $C>0$ such that for all $\lambda>1$,
any mapping $p\mapsto q = q(p)$ satisfying $q(p) > \max(\lambda p,p+4)$
yields
\begin{equation}\label{saturation-3}
\Vert u \Vert_{H^1_0(\Omega)} \leq 
C \frac{\lambda}{\lambda-1} \, \Vert u_q \Vert_{H^1_0(\Omega)} \qquad \textrm{for all } \, p\geq 0 \text{ and } 
 f \in {\mathbb P}_p(\Omega).
\end{equation}
\end{theorem}

Since most $hp$-FEMs in the literature perform $p$-enrichment upon
adding a constant increment to $p$, typically $1$ or $2$, one may
wonder whether the preceding sufficient condition on $q$ is also
necessary. We now investigate this question computationally upon
defining
\begin{equation*}
  C_{p,q,r} :=\max_{f\in\rhn{\mathbb{P}_p(\Omega)}}
  \frac{\|u_r\|_{H^1_0(\Omega)}}{\|u_q\|_{H^1_0(\Omega)}}
\end{equation*}
where $r\gg q$ is chosen computationally so that $u_r$ is sufficiently close to $u$ in
$H^1_0(\Omega)$; note that this is not a hidden saturation assumption
because the value of $r$ is not predetermined but found once
the number $C_{p,q,r}$ stabilizes. This calculation reduces to an eigenvalue problem,
already used in \cite{CNSV17}, and leads to Figure
\ref{F:constant-increment} for $q=p+k$ with $k=2,4,6,10$:
\begin{figure}[h!]
\begin{center}
\includegraphics[width=9.0cm]{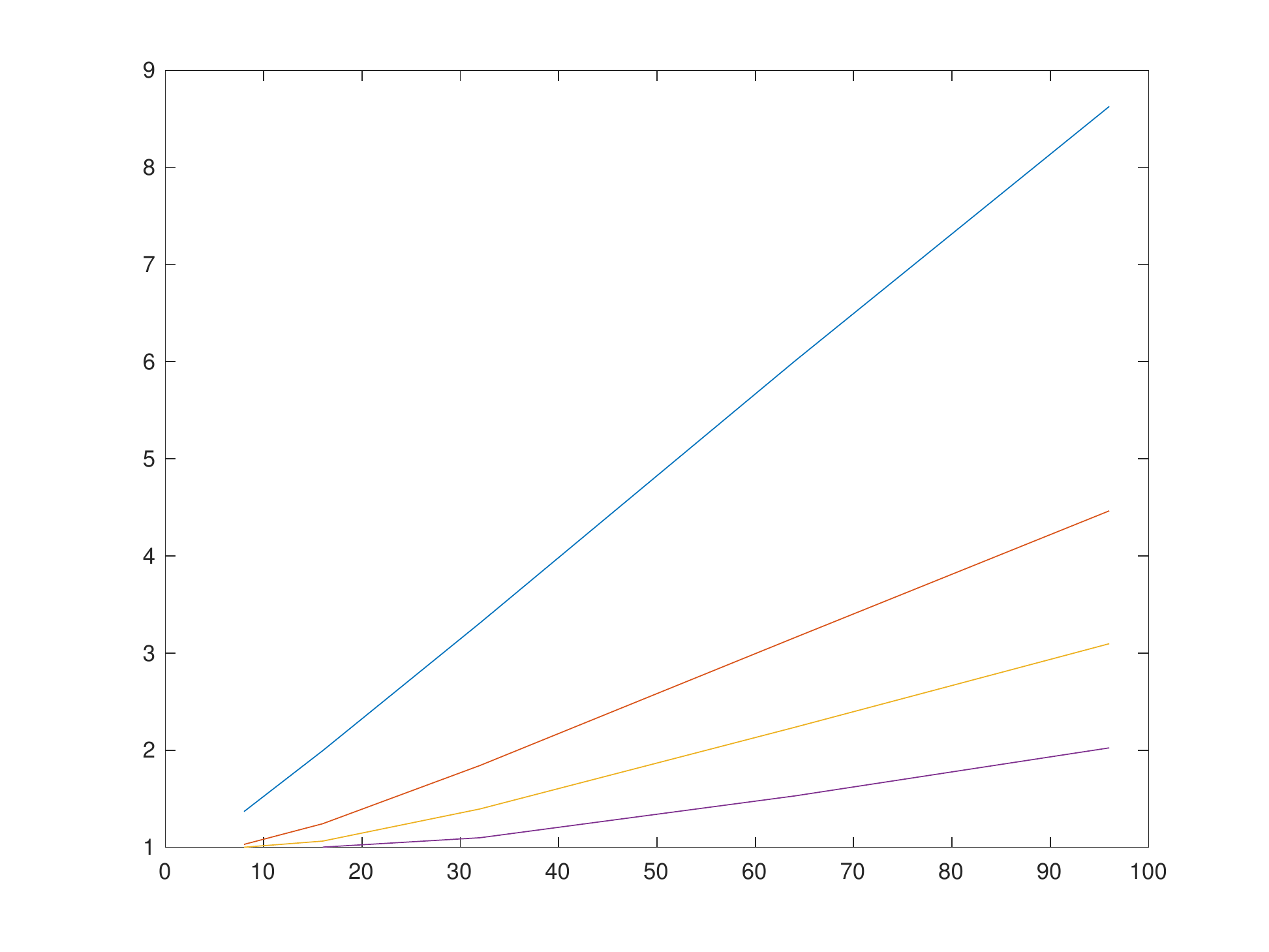}
\caption{Constants $C_{p, q, r}$ vs \rhn{$q=p+k$: blue $k=2$, orange $k=4$,
  yellow $k=6$, violet $k=10$.} The dependence is clearly linear
  rather than constant but the growth is moderate.}
\label{F:constant-increment}
\end{center}
\end{figure}
We thus realize that $C_{p,q,r}$ exhibits a modest but linear growth
on $q=p+k$ for $k$ constant, which confirms that this choice is not
$p$-robust. For moderate values of $p$ this might still be acceptable
computationally, but it could compromise computational complexity for
extreme values of $p$ as in spectral algorithms \cite{CHQZ06}.

Even though the saturation property is quite delicate, it has been
often used in a posteriori error analysis of low order AFEMs until
now. It originates in the work of Bank and Weiser \cite{BW85}, and
Bornemann, Erdmann and Kornhuber \cite{BEK:96}; see Nochetto
\cite{Noch93} for related work. D\"orfler and Nochetto \cite{DN02} proved
the saturation property for $p=1, q=2$ and $d=2$ provided data
oscillation is small relative to $\|u-u_p\|_{H^1_0(\Omega)}$ but
showed counterexamples for piecewise constant forcing $f$.

We stress that \eqref{saturation-3} is not asymptotic: it is valid for
any $p\ge0$ and any $f\in\mathbb{P}_p(\Omega)$. Since $u_q\to u$ in
$H^1_0(\Omega)$ as $q\to\infty$ it is obvious that $C_{p,q,\infty}\to
1$ as $q\to\infty$. It is for this reason that Theorem
\ref{theo:main} has some intrinsic value in the theory of FEMs and
might have implications beyond a posteriori error analysis.

The proof of Theorem \ref{theo:main} proceeds as follows. We perform a
multilevel decomposition of $\cV_q$
\begin{equation}\label{multilevel-decom}
\cV_q = \bigoplus_{j=1}^q \cW_j,
\end{equation}
where $\cW_j$ are polynomial subspaces of total degree $j$. Since this
decomposition is quasi-orthogonal in the sense that
\[
\cW_j \perp \cW_\ell
\quad\textrm{for all } \ell\ne j-2, j, j+2,
\]
we need to account for interactions between neighboring spaces
$\cW_j$. We study the angle between subspaces $\cW_j$ and show it is
larger than $\pi/3$; this is the content of Proposition
\ref{prop:proj}. This in turn allows us to find the precise decay of high
frequency modes of $u_q$, which leads to \eqref{saturation-2}.

The paper is organized as follows. In section \ref{sec:bases} we
introduce the multilevel decomposition \eqref{multilevel-decom} and
discuss a few properties including Proposition \ref{prop:proj}.
In section \ref{sec:decay} we analyze
the decay of high order components of $u_q$, whereas in section
\ref{sec:proof} we prove Theorem \ref{theo:main}.
We conclude in section \ref{App} with the proof of
the rather technical Proposition \ref{prop:proj}.


\section{Multi-level decompositions of polynomial spaces}\label{sec:bases}

Hereafter, we recall the definition of classical polynomial bases in $L^2(\Omega)$ and in $H^1_0(\Omega)$, obtained by tensorization from
corresponding bases in $L^2(I)$ and in $H^1_0(I)$, where $I=(-1,1)$ is the reference interval. The elements of theses bases enjoy certain orthogonality properties, by which a multi-level, quasi-orthogonal decomposition of $H^1_0(\Omega)$ is obtained. This will be useful in deriving the main result of this paper.

On the interval $I$, we consider the {\sl orthonormal Legendre basis in $L^2(I)$}
\begin{equation}\label{eq:A1}
\vartheta_k(x) =\sqrt{k+\tfrac12} \, L_k(x), \qquad k \geq 0,
\end{equation}
(where $L_k$ stands for the $k$-th Legendre orthogonal polynomial in $I$, 
which satisfies  ${\rm deg}\, L_k = k$ and $L_k(1)=1$), as well as the
{\sl orthonormal Babu\v ska-Shen \rhn{(BS)} basis in $H^1_0(I)$}:
\begin{equation}\label{eq:A2}
\begin{split}
\varphi_k({x})&=\sqrt{k-\tfrac12}\int_{{x}}^1 L_{k-1}(s)\,{d}s \\
&=\frac1{\sqrt{4k-2}}\big(L_{k-2}({x})-L_{k}({x})\big), 
 \qquad k \geq 2\;.
 \end{split}
\end{equation}
The BS basis enjoys the following orthogonality properties in $L^2(I)$
\rhn{for $m\ge k$}:
\begin{equation}\label{eq:A3}
(\varphi_k,\varphi_m)_{L^2({I})} = 
\begin{cases}
\frac2{(2k-3)(2k+1)} & \text{if } m=k \;, \\
- \frac1{(2k+1)\sqrt{(2k-1)(2k+3)}} &  \text{if } m=k+2 \;, \\
0 & \text{otherwise.}
\end{cases}
\end{equation}

On the square $\Omega=I \times I$, the previous bases induce, resp., the {\sl tensorized orthonormal Legendre basis in $L^2(\Omega)$}:
\begin{equation}\label{eq:A4}
\Theta_k(x)= \vartheta_{k_1}(x_1)\vartheta_{k_2}(x_2), \qquad k \in \hat{\cal K},
\end{equation}
where $k=(k_1,k_2)$, $x=(x_1, x_2)$ and  $\hat{\cal K}=\mathbb{N}^2$, and the {\sl tensorized Babu\v ska-Shen basis in  $H^1_0(\Omega)$}:
\begin{equation}\label{eq:A5}
\Phi_k(x) = \varphi_{k_1}(x_1) \varphi_{k_2}(x_2),  \qquad k \in {\cal K}\;, 
\end{equation}
where ${\cal K}=\{k \in \mathbb{N}^2 \ : \ k_i \geq 2 \text{ for }i=1,2 \}$.  

The tensorized BS basis is not orthogonal in $H^1_0(\Omega)$. Indeed, from the expression
\rhn{
\begin{align*}
(\Phi_k,\Phi_m)_{H^1_0(\Omega)} &=
(\varphi_{k_1},\varphi_{m_1})_{H^1_0({I})}(\varphi_{k_2},\varphi_{m_2})_{L^2({I})}
\\ &+
 (\varphi_{k_1},\varphi_{m_1})_{L^2({I})}(\varphi_{k_2},\varphi_{m_2})_{H^1_0({I})} \;,
\end{align*}
}
and (\ref{eq:A3}) we immediately obtain 
\begin{equation}\label{eq:A6}
(\Phi_k,\Phi_m)_{H^1_0(\Omega)} \not = 0 \text{ \ iff \ }
\begin{cases} 
k_1=m_1 \text{ and } k_2-m_2 \in \{-2,0,2\}, \text{ or } \\
k_2=m_2 \text{ and } k_1-m_1 \in \{-2,0,2\}.
\end{cases} 
\end{equation}
As a consequence, denoting by $|k|=|k_1|+|k_2|$ the $\ell^1$-norm in $\mathbb{Z}^2$,  we have
\begin{equation}\label{eq:A7}
(\Phi_k,\Phi_m)_{H^1_0(\Omega)} = 0 \qquad \text{if  \ } |k|-|m| \not\in \{-2,0,2\}.
\end{equation}

At last, concerning the interaction between  the Legendre basis and the BS one, we have
\begin{equation}\label{eq:A8} 
(\Theta_k,\Phi_m)_{L^2(\Omega)} \not = 0 \text{ \ iff \ } k_1 \in \{m_1-2,m_1\} \text{ and } k_2 \in \{m_2-2,m_2\},
\end{equation}
which implies
\begin{equation}\label{eq:A9} 
(\Theta_k,\Phi_m)_{L^2(\Omega)} = 0 \qquad \text{if  \ } |k-m| > 4.
\end{equation}

\medskip
\begin{remark}[{\sl orthogonality by parity}]\label{rem:parity}
Any function $v \in L^2(\Omega)$ can be \rhn{split uniquely} into four components
\begin{equation}\label{eq:parity}
v = \sum_{\alpha \in \{0,1\}^2} v^\alpha,
\end{equation}
where $v^\alpha$ for $\alpha=(\alpha_1, \alpha_2)$ is even (odd, resp.) with respect to the variable $x_i$ ($i=1,2$) iff $\alpha_i=0$ ($\alpha_i=1$, resp.). 
For instance, $v^{(0,1)}$ satisfies $v^{(0,1)}(-x_1,x_2)=v^{(0,1)}(x_1,x_2)$ and $v^{(0,1)}(x_1,-x_2)=-v^{(0,1)}(x_1,x_2)$ for all
$(x_1,x_2)\in \Omega$.

Components with different parity indices are always $L^2(\Omega)$-orthogonal, and $H^1_0(\Omega)$-orthogonal whenever $v \in H^1_0(\Omega)$. In particular, as a consequence of \eqref{eq:A6} and \eqref{eq:A8}, we observe that $(\Phi_k,\Phi_m)_{H^1_0(\Omega)} = 0$ and $(\Theta_k,\Phi_m)_{L^2(\Omega)} = 0$ whenever $k$ and $m$ have at least one entry of different parity.
\end{remark}

\subsection{Detail spaces and their projectors}

For $j \geq 4$, let us define the finite dimensional subspace of $H^1_0(\Omega)$
\begin{equation}\label{eq:defWj}
\cW_j := \text{span}\{\Phi_k \, :  |k|=j \}\;.
\end{equation}
Note that, thanks to (\ref{eq:A6}), the functions $\Phi_k$ that generate $\cW_j$ are mutually orthogonal in $H^1_0(\Omega)$. 
We immediately have the multi-level decompositions
\begin{equation}\label{eq:multilevel}
\cV_q = \bigoplus_{j=4}^q \cW_j  \quad \textrm{for all } q \geq 4  \qquad \text{and} \qquad H^1_0(\Omega) = \bigoplus_{j=4}^\infty \cW_j \;.
\end{equation}
Such decompositions are `quasi-orthogonal', in the sense that by (\ref{eq:A7}) we have
\begin{equation} \label{eq:Wj-orth}
\cW_j \perp_{H^1_0(\Omega)} \cW_\ell \qquad \text{for all \ } \ell  \not = j-2, j, j+2.
\end{equation}
Furthermore, the `angle' between two non-orthogonal subspaces is uniformly bounded away from 0, as implied by the
following \rhn{technical result, that will be crucial in the sequel.
We postpone its proof to section \ref{App}.}
\begin{proposition}[angle between $\cW_{j-2}$ and $\cW_j$]\label{prop:proj} 
Let $P_j : \cW_{j-2} \to \cW_j$ ($j \geq 6$) be the orthogonal projection with respect to the $H^1_0(\Omega)$-inner product. 
Then,
$$
\Vert P_j \Vert_{{\cal L}(\cW_{j-2},\cW_j)} < \frac12.
$$
Actually, there exists a constant $c>0$ independent of $j$ such that 
$$
\Vert P_j \Vert_{{\cal L}(\cW_{j-2},\cW_j)} \leq \frac12 \left( 1 - \frac{c}{\, j^{2}} \right) \;. 
$$
\end{proposition}

Note that the orthogonal projection $P_j^*: \cW_{j} \to \cW_{j-2}$, given by the adjoint of $P_j$, satisfies the same estimate.

\section{Decay of the higher-order components of the Galerkin solution}\label{sec:decay}

Given $f \in {\mathbb P}_p(\Omega)$, let  $u_q \in \cV_q$ be the
Galerkin solution defined in \eqref{eq:def-uq}, and let $u_q =
\sum_{j=4}^q U_j$, with $U_j=U_j(q) \in \cW_j$, be its multilevel decomposition according to \eqref{eq:multilevel}. The purpose of this section is to prove that for any $q$ sufficiently larger than $p$, the $H^1_0(\Omega)$-norm of  $U_q$ and $U_{q-1}$ decay at least proportionally to the quantity $(q-p)^{-1}$. The precise result is as follows.

\begin{proposition}[decay of $U_j$]\label{prop:p1}
For any $p \geq 0$ and $q > \hat p := p+4$, one has 
\begin{equation}\label{eq:c1}
\Vert U_j \Vert_{H^1_0(\Omega)}  \leq \frac6{q-p} \, \Vert u_q \Vert_{H^1_0(\Omega)}\;,  \qquad j=q, \, q-1.
\end{equation}
\end{proposition}

\proof
We first observe that the parity splitting \eqref{eq:parity} of the forcing $f$ induces by linearity a corresponding splitting of the Galerkin solution $u_q$ as well as of each of its multi-level details $U_j$, which is nothing but the parity splitting of $u_q$ as well as of $U_j$. Therefore, thanks  to the orthogonality of the components with different parity (cf. Remark \ref{rem:parity}), it is enough to establish \eqref{eq:c1} for each component separately, and then sum-up the squares of both sides invoking Parseval's identity.

For the sake of definiteness, we will focus on the components of (even, even) type, the other types being amenable to a similar treatment. Thus, referring to \eqref{eq:parity} for the notation, we consider the component $u_q^{(0,0)}$ of $u_q$ (which solves \eqref{eq:def-uq} for the forcing $f^{(0,0)}$), as well as its details $U_j^{(0,0)} \in \cW_j^{(0,0)}:= \text{span}\{\Phi_k \, :  |k|=j \text{ and } k_1, k_2 \text{ are even}\}$. We aim at proving that for $q \geq \hat p+1$ and $j \in \{q-1,q\}$
$$
\Vert U_j^{(0,0)} \Vert_{H^1_0(\Omega)} \leq \frac6{q-p} \, \Vert u_{q}^{(0,0)} \Vert_{H^1_0(\Omega)}.
$$
However, it is easily seen that $U_{q-1}^{(0,0)}=0$ if $q$ is even, and similarly  $U_{q}^{(0,0)}=0$ if $q$ is odd. Hence, we will prove
\begin{equation}\label{eq:c11}
\Vert U_q^{(0,0)} \Vert_{H^1_0(\Omega)} \leq \frac6{q-p} \, \Vert u_{q}^{(0,0)} \Vert_{H^1_0(\Omega)}
\end{equation}
under the assumption that $q$ is even, the other situation being similar.

To avoid cumbersome notation, for the rest of the proof we will drop the superscript $\phantom{u}^{(0,0)}$ from all entities. 
So, we will write
\[
u_q=\sideset{}{'}\sum\limits_{j=4}^q U_j
\qquad\textrm{with } U_j \in \cW_j,
\]
where here and in the sequel the symbol ${}^{'}$ indicates that the summation runs over even indices only.

From the Galerkin equations, we have for any even $j \in [4,q]$
\begin{equation} \label{eq:r1}
(u_q,W_j)_{H^1_0(\Omega)}=(f,W_j)_{L^2(\Omega)} \qquad \textrm{for all
  } W_j\in \cW_j \; .
\end{equation}
Since $q \geq \hat{p}+1$, exploiting (\ref{eq:A7}) and (\ref{eq:A8}),
\eqref{eq:r1} yields
\begin{equation}
  (U_q,W_q)_{H^1_0(\Omega)}+(U_{q-2},W_q)_{H^1_0(\Omega)}=0\qquad
 \textrm{for all } W_q\in \cW_q \;, 
\end{equation}
which can be rewritten equivalently as 
\begin{equation}\label{eq:j=q}
U_q=-T_q^{-1}P_q U_{q-2}
\end{equation}
where $T_q=I$ and $P_q$ is defined in Proposition \ref{prop:proj}.

For any even $j$ satisfying $\hat{p}+2<j \leq q-2$, \eqref{eq:r1} yields
\begin{equation} \label{eq:r2}
  (U_{j+2},W_j)_{H^1_0(\Omega)}+(U_{j},W_j)_{H^1_0(\Omega)}+(U_{j-2},W_j)_{H^1_0(\Omega)}=0
  \qquad 
\end{equation}
for all $W_j\in \cW_j$;
this is equivalent to  $P_{j+2}^* U_{j+2} + U_j + P_j U_{j-2} = 0$.
Assuming by induction that $U_{j+2}=-T_{j+2}^{-1} P_{j+2} U_j$ with
$\|T_{j+2}^{-1}\| \leq 2$ (which is trivially true for $j=q-2$
according to \eqref{eq:j=q}), gives 
$$
(I-P_{j+2}^* T_{j+2}^{-1} P_{j+2}) U_j=-P_j U_{j-2}.
$$ 
Using Proposition~\ref{prop:proj}, the operator
\begin{equation}\label{def:Tj}
T_j:=I-P_{j+2}^* T_{j+2}^{-1} P_{j+2}
\end{equation}
is invertible and satisfies
\begin{equation}\label{eq:norminv}
 \| T^{-1}_j\| \leq \frac {1}{1- \frac 1 4   \| T^{-1}_{j+2}\| }\;.
 \end{equation}
We conclude that
\begin{equation} \label{eq:r3}
U_j=-T_j^{-1}P_j U_{j-2}.
\end{equation}
with $\| T^{-1}_{j}\| \leq 2$, which proves the induction argument for all even $j$ satisfying $\hat{p} + 2 < j \le q-2$.

Next, we have to bound the norm of $U_{j-2}$ for $j=\hat{p}+4$ when $p$, hence $\hat{p}$, is even, or for $j=\hat{p}+3$ when $p$ is odd. It is therefore convenient to define the even integer
$$
r:=\begin{cases} \hat{p} = p+4 & \text{if } p \text{ is even}, \\
\hat{p}-1 = p+3 & \text{if } p \text{ is odd}, 
\end{cases}
$$
so that in both cases, we have to bound $\|U_{r+2}\|_{H^1_0(\Omega)}$.
To this end, let us introduce
\[
\cV_r^{(0,0)}:= \sideset{}{'}\bigoplus_{j=4}^r \cW_j,
\qquad
\bar u_r:=\sideset{}{'}\sum\limits_{j=4}^r
  U_j\in\cV_r^{(0,0)};
\]
note that $\bar u_r \ne u_r$ because $U_j=U_j(q)$.
Then, in view of \eqref{eq:r1}, we deduce
\begin{equation}
(\bar u_r+\sideset{}{'}\sum_{j=r + 2}^q U_j, v_r)_{H^1_0(\Omega)}=(f,v_r)_{L^2(\Omega)} \quad \textrm{for
    all } v_r\in \cV_r^{(0,0)} \;,
\end{equation}
which, thanks to $\displaystyle{(\sideset{}{'}\sum_{j=r + 2}^q U_j, v_r)_{H^1_0(\Omega)}}=(P_r^*U_{r+2},v_r)_{H^1_0(\Omega)}$ for all $v_r \in \cV_r^{(0,0)}$, implies that
$$
\bar u_r+P_r^* U_{r+2}=u_r.
$$

We observe that \eqref{eq:r2} is also valid for $j=r+2$.
Since $(U_r,W_{r+2})_{H^1_0(\Omega)}= (\bar u_r,W_{r+2})_{H^1_0(\Omega)}$
we obtain
\begin{equation}
(U_{r +4},W_{r +2})_{H^1_0(\Omega)}+(U_{r +2},W_{r +2})_{H^1_0(\Omega)}+(\bar u_r,W_{r +2})_{H^1_0(\Omega)}= 0
\end{equation}
and
$$
T_{r +2} U_{r +2} = -P_{r +2} \bar u_r 
$$
as it happened with \eqref{eq:r3}. This implies 
\begin{equation}
(T_{r +2}-P_{r +2}P_r^*) U_{r +2}=-P_{r +2} u_r
\end{equation}
which in view of \eqref{def:Tj} yields
\begin{equation}
(I-P_{r +4}^* T_{r +4}^{-1} P_{r +4} -P_{r +2}P_r^*) U_{r +2}=-P_{r +2} u_r.
\end{equation}
Since $\|P_{r +4}^* T_{r +4}^{-1} P_{r +4} +P_{r +2}P_r^*\|\leq 2
  \frac{1}{4}+\frac{1}{4}=\frac{3}{4}$, thanks to Proposition~\ref{prop:proj} and $\|T_{r +4}^{-1}\|\le2$, we conclude that
$$
\|U_{r +2}\|_{H^1_0(\Omega)} \leq \frac{1}{1-\frac{3}{4}}\frac{1}{2} \|u_r\|_{H^1_0(\Omega)}=2\|u_r\|_{H^1_0(\Omega)} \leq 2 \Vert u_q\Vert_{H^1_0(\Omega)},
$$
where last inequality follows from the inclusion
$\cV_r^{0,0)} \subset \cV_q^{(0,0)}$ and the
minimization property of the Galerkin solution.
%
%

Collecting the above results we arrive at
\begin{equation}\label{eq:collect}
\begin{split}
& U_j = -T_j^{-1}P_j U_{j-2} \qquad r +4 \leq   j \leq q  \quad \text{($j$ even)}, \\
& \| U_{r +2}\|_{H^1_0(\Omega)} \leq 2 \Vert u_q \Vert_{H^1_0(\Omega)}.
\end{split} 
\end{equation}
For $q \geq \hat{p}+4$, this implies 
$$
\|U_q\|_{H^1_0(\Omega)} \leq
\sideset{}{'}\prod_{j=r +4}^{q} \| T_j^{-1}P_j \| \| U_{r +2}\|_{H^1_0(\Omega)}
\leq   2  \Vert u_q \Vert_{H^1_0(\Omega)} \sideset{}{'}\prod_{j=r+4}^{q} \| T_j^{-1}P_j \|.
$$
In order to bound the product on the right-hand side, let us write $j=q-2m$ with
$m=0,1, \dots, s$ and $s:= \frac12(q-r)-2$.
Then, by Proposition \ref{prop:proj}, we have $\| T_j^{-1}P_j \| \leq \frac12 \| T_j^{-1}\| =: \alpha_m$. Recalling \eqref{eq:norminv}, it holds
$$
\alpha_m \leq \frac{\frac12}{1-\frac12 \alpha_{m-1}}= \frac1{2-\alpha_{m-1}}, \qquad \text{with \ } \alpha_0 \leq \tfrac12 \;.
$$
By recurrence, it is immediate to check that $\alpha_m \leq \frac{m+1}{m+2}$, whence
$$
 \sideset{}{'}\prod_{j=r+4}^{q} \| T_j^{-1}P_j \| \leq  \prod_{m=0}^s \alpha_m \leq \prod_{m=0}^s \frac{m+1}{m+2} =\frac1{s+2}=
\frac2{q-r}.
$$
Since $q\ge p+6$ if $p$ is even and $q\ge p+5$ if $p$ is odd,
it is easily checked that
$$
\frac2{q-r} =\begin{cases} \frac2{(q-p)-4} \leq \frac6{q-p} & \text{if } p \text{ is even}, \\
\frac2{(q-p)-3} \leq \frac5{q-p} & \text{if } p \text{ is odd}. 
\end{cases}
$$
This gives the desired estimate \eqref{eq:c11}.
\endproof

\section{A subspace decomposition in $H^1_0(\Omega)$  }\label{sec:angle}

Consider the complementary space of $\cV_q$ in $H^1_0(\Omega)$ given by
\begin{equation}\label{eq:compl-space}
\cV_q^c :=\clos_{H^1_0(\Omega)} \text{span}\, \{ \Phi_m \, : \, |m| > q \}.
\end{equation}
Therefore, $H^1_0(\Omega)=\cV_q \oplus \cV_q^c$ and 
any $v \in H^1_0(\Omega)$ can be split as
\[
v = v_q + z_q, \qquad
v_q \in \cV_q, \quad z_q \in \cV_q^c.
\]
The purpose of this section is to apply once more Proposition \ref{prop:proj} and derive a bound on the norm of $v_q$ and $z_q$ in terms of the norm of $v$. 

We start with the following auxiliary result for any $w=\sum_{j=4}^q W_j\in\cV_q.$

\begin{lemma}[bound of $\Vert W_q \Vert_{H^1_0(\Omega)}$]\label{lem:proj}
For any $q \geq 4$ and any $w=\sum_{j=4}^q W_j\in\cV_q$, one has
$$
\Vert W_q \Vert_{H^1_0(\Omega)} \leq \sqrt{2} \Vert w \Vert_{H^1_0(\Omega)} \;.
$$ 
\end{lemma}
\proof As in the previous section, splitting $w$ and $W_q$ in their orthogonal components according to the parity of the basis functions, it is enough to establish the result for each component separately. Hereafter, we detail the analysis for the `(even, even)' case, in which case we may assume $q$ even, since otherwise $W_q^{(0,0)}=0$ and the result is trivial. 

Dropping as above the superscript $\phantom{u}^{(0,0)}$ in functions
and subspaces, we write $w=W+W_q$ with
\[
W = \sideset{}{'}\sum\limits_{j=4}^{q-2} W_j \in \cV_{q-2}.
\]
Keeping $W_q$ fixed, let us first minimize the norm of $w$, i.e., let
us look for the minimizer $\bar W\in \cV_{q-2}$ of the quantity
$\Psi(W):=\Vert W+W_q\Vert_{H^1_0(\Omega)}^2$. Such a function
satisfies
\begin{equation}\label{eq:staz1}
(\bar W, Y)_{H^1_0(\Omega)} = - (W_q, Y)_{H^1_0(\Omega)} \qquad
  \textrm{for all } Y \in \cV_{q-2}
\end{equation}
and
\begin{equation}\label{eq:staz2}
\Psi(\bar W) = \Vert W_q \Vert_{H^1_0(\Omega)}^2 +(W_q, \bar W)_{H^1_0(\Omega)} \;.
\end{equation}
Using the orthogonality conditions \eqref{eq:A7}, we obtain the sequence of equations
$$
(\bar W_4,Y_4)_{H^1_0(\Omega)}+(\bar W_6,Y_4)_{H^1_0(\Omega)}=0\qquad
\textrm{for all } Y_4\in \cW_4\;,
$$
and
$$
(\bar W_{j-2},Y_j)_{H^1_0(\Omega)}+(\bar
W_{j},Y_j)_{H^1_0(\Omega)}+(\bar W_{j+2},Y_j)_{H^1_0(\Omega)}=0\qquad
\textrm{for all } Y_j\in \cW_j
$$
for any even $j$ such that $4 < j <q-2$, and finally
$$
(\bar W_{q-4},Y_{q-2})_{H^1_0(\Omega)}+(\bar
W_{q-2},Y_{q-2})_{H^1_0(\Omega)}=-( W_q,Y_{q-2})_{H^1_0(\Omega)}
$$
for all $Y_{q-2}\in \cW_{q-2}$.
Setting recursively $T_4=I$ and $T_j=(I-P_j T_{j-2}^{-1} P_j^*)$,  we derive $\bar W_j =-T_j^{-1} P_{j+2}^* \bar W_{j+2}$ for $j=4,6, \dots, q-4$, and $\bar W_{q-2}=-T_{q-2}^{-1} P_{q}^* W_q$. Note that, thanks to Proposition \ref{prop:proj}, one can prove as in Sect. \ref{sec:decay} that $\Vert T_j^{-1} \Vert \leq 2$ for all $j$. Since
\begin{align*}
(W_q, \bar W)_{H^1_0(\Omega)} &= (W_q, \bar W_{q-2})_{H^1_0(\Omega)} \\
  &=  (P_q^* W_q, \bar W_{q-2})_{H^1_0(\Omega)}
  = - (P_q^* W_q, T_{q-2}^{-1} P_{q}^* W_{q})_{H^1_0(\Omega)},
\end{align*}
using once more Proposition \ref{prop:proj}, we deduce
\begin{align*}
(P_q^* W_q, T_{q-2}^{-1} P_{q}^* W_{q})_{H^1_0(\Omega)} &\leq \Vert T_{q-2}^{-1} \Vert \, \Vert P_q^* \Vert^2 \Vert W_q \Vert_{H^1_0(\Omega)}^2 \\
&\leq 2 \left(\frac12\right)^2 \Vert W_q \Vert_{H^1_0(\Omega)}^2 =\frac12 \Vert W_q\Vert_{H^1_0(\Omega)}^2\;.
\end{align*}
In view of \eqref{eq:staz2} and the preceding estimate,
we conclude that 
$$
\Vert w \Vert_{H^1_0(\Omega)}^2 = \Psi(W) \geq \Psi(\bar W) \geq \Vert W_q \Vert_{H^1_0(\Omega)}^2 -\frac12 \Vert W_q \Vert_{H^1_0(\Omega)}^2 = \frac12 \Vert W_q \Vert_{H^1_0(\Omega)}^2
$$
for any $w \in \cV_q^{(0,0)}$, whence the asserted estimate follows.
\endproof

We now establish the main result of this section.

\begin{proposition}[control of $\Vert z_q \Vert_{H^1_0(\Omega)}$]\label{prop:p2}
There exists a constant $C_2>0$ such that for any $q \geq 4$ and any
$v=v_q + z_q \in \cV_q \oplus \cV_q^c$, one has 
\begin{equation}\label{eq:c2}
\Vert z_q \Vert_{H^1_0(\Omega)} \leq C_2  \, q \, \Vert v \Vert_{H^1_0(\Omega)}.
\end{equation}
\end{proposition}

\proof Using $\Vert z_q \Vert_{H^1_0(\Omega)} \leq \Vert v \Vert_{H^1_0(\Omega)} + \Vert v_q \Vert_{H^1_0(\Omega)}$, it is enough to prove the existence of a constant $C_2'>0$ independent of $q$ such that for all $v \in H^1_0(\Omega)$ 
\begin{equation}\label{eq.c2primo}
\Vert v_q \Vert_{H^1_0(\Omega)} \leq C_2' q \Vert v \Vert_{H^1_0(\Omega)}\;.
\end{equation}
To this end, let us focus as above on the `(even, even)' components of
$v$ and $v_q$, in which case it is not restrictive to assume $q$ even,
and drop the superscript $\phantom{u}^{(0,0)}$ in functions and subspaces.
Let us fix any even integer $r>q$ and assume first that $v \in \cV_r$
is written as $v=v_q+V$, with
\[
V= \sideset{}{'}\sum\limits_{j=q+2}^r V_j \in
\sideset{}{'}\bigoplus\limits_{j=q+2}^r \cW_j.
\]
By applying the same technique as above, i.e., minimizing the (squared) norm $\Psi(V):=\Vert v_q+V\Vert_{H^1_0(\Omega)}^2$, we find that
\begin{equation}\label{eq:lwbd}
\Vert v \Vert_{H^1_0(\Omega)}^2 = \Psi(V) \geq \Psi(\bar V) = \Vert v_q \Vert_{H^1_0(\Omega)}^2 +(v_q, \bar V)_{H^1_0(\Omega)} \;,
\end{equation}
where the minimizer $\bar V =  \sideset{}{'}\sum\limits_{j=q+2}^r \bar V_j$ is such that $\bar V_{q+2}=-T_{q+2}^{-1} \tilde P_{q+2} v_q$, for 
$T_{q+2}$ defined recursively by \eqref{def:Tj} with $T_r=I$, and $\tilde P_{q+2} : H^1_0(\Omega) \to \cW_{q+2}$ defined as the orthogonal projection in the $H^1_0(\Omega)$ inner product. Now,
\begin{align*}
(v_q, \bar V)_{H^1_0(\Omega)} &= (v_q, \bar V_{q+2})_{H^1_0(\Omega)} = (\tilde P_{q+2}v_q, \bar V_{q+2})_{H^1_0(\Omega)} \\
&= - (\tilde P_{q+2}v_q, T_{q+2}^{-1} \tilde P_{q+2} v_q)_{H^1_0(\Omega)}
\end{align*}
with
$$
\big| (\tilde P_{q+2}v_q, T_{q+2}^{-1} \tilde P_{q+2} v_q)_{H^1_0(\Omega)} \big| \leq \Vert T_{q+2}^{-1} \Vert \, \Vert  \tilde P_{q+2} v_q \Vert_{H^1_0(\Omega)}^2 \leq 2 \Vert  \tilde P_{q+2} v_q \Vert_{H^1_0(\Omega)}^2\;.
$$ 
Writing $v_q =\sideset{}{'}\sum\limits_{j=4}^q V_j$, one has $\tilde
P_{q+2} v_q = \tilde P_{q+2} V_q =  P_{q+2} V_q $, whence by
Proposition \ref{prop:proj} with $\varepsilon_j=cj^{-2}$
and Lemma \ref{lem:proj} we get
\begin{align*}
\Vert  \tilde P_{q+2} v_q \Vert_{H^1_0(\Omega)} = \Vert  P_{q+2} V_q \Vert_{H^1_0(\Omega)} 
&\leq \frac12 \left(1-\varepsilon_{q+2}\right) \Vert  V_q \Vert_{H^1_0(\Omega)} \\ &\leq \frac1{\sqrt{2}} \left(1-\varepsilon_{q+2}\right) \Vert  v_q \Vert_{H^1_0(\Omega)} \;,
\end{align*}
which gives
$$
\big| (\tilde P_{q+2}v_q, T_{q+2}^{-1} \tilde P_{q+2} v_q)_{H^1_0(\Omega)} \big| \leq \left(1-\varepsilon_{q+2}\right)^2 \Vert  v_q \Vert_{H^1_0(\Omega)}^2 \;.
$$
Then, from \eqref{eq:lwbd} we obtain
$$
\Vert v \Vert_{H^1_0(\Omega)}^2 \geq \varepsilon_{q+2}(2-\varepsilon_{q+2}) \Vert  v_q \Vert_{H^1_0(\Omega)}^2\;,
$$
which immediately yields \eqref{eq.c2primo} for all
$v \in \cV_r=\cV_r^{(0,0)}$ and all $r>q$. 

The same result holds for all other combinations of parity indices; hence, it holds for any $v \in \cV_r$. Since polynomials vanishing on $\partial\Omega$ form a dense subset of $H^1_0(\Omega)$, we conclude that \eqref{eq.c2primo} holds for all $v \in H^1_0(\Omega)$. \endproof

\section{Proof of Theorem \ref{theo:main}}\label{sec:proof}

We actually prove the equivalent condition
$$
\Vert u - u_q \Vert_{H^1_0(\Omega)} \lesssim  \Vert u_q \Vert_{H^1_0(\Omega)},
$$
and for that we write
$$
\Vert u - u_q \Vert_{H^1_0(\Omega)}  = \sup_{v \in H^1_0(\Omega), \ v \not = 0} 
\frac{ (u-u_q, v)_{H^1_0(\Omega)}}{\Vert  v \Vert_{H^1_0(\Omega)}}.
$$
As in the previous section, let us split any $v \in H^1_0(\Omega)$ as $v = v_q + z_q \in \cV_q \oplus \cV_q^c$, where $\cV_q^c$ is given by (\ref{eq:compl-space}). By the Galerkin orthogonality and the definition of $u$, we have
$$
(u-u_q, v)_{H^1_0(\Omega)} = (u-u_q, z_q)_{H^1_0(\Omega)} = (f,z_q)_{L^2(\Omega)} - (u_q, z_q)_{H^1_0(\Omega)} .
$$
Recalling (\ref{eq:A8}) and the condition $q > \hat{p}$, we have $ (f,z_q)_{L^2(\Omega)} =0$, hence
$$
(u-u_q, v)_{H^1_0(\Omega)} =  - (u_q, z_q)_{H^1_0(\Omega)} .
$$
Now, recalling (\ref{eq:multilevel}), we expand the Galerkin solution $u_q$ as $u_q =\sum_{j=4}^q U_j$. Invoking (\ref{eq:A7}), we get
$$
(u_q, z_q)_{H^1_0(\Omega)} = (U_{q-1}+U_q, z_q)_{H^1_0(\Omega)}. 
$$
Applying Propositions \ref{prop:p1} and \ref{prop:p2}, we get the following bound
$$
(u-u_q, v)_{H^1_0(\Omega)} \leq C_1C_2 \frac{q}{q-p} \Vert u_q \Vert_{H^1_0(\Omega)}  \Vert v \Vert_{H^1_0(\Omega)}.
$$
Since for $q >\lambda p$, \rhn{the relation $\frac{q}{q-p}
<\frac{\lambda}{\lambda-1}$ holds,
and the proof is complete.} \endproof

\section{Proof of Proposition \ref{prop:proj}}\label{App}

We establish \rhn{the bound $\Vert P_{j+2}\Vert_{{\cal L}(\cW_j,\cW_{j+2})}
\le \frac{1}{2} \big(\frac{1}{2} - \frac{c}{j^2}\big)$
for any $j \geq 4$. This will be achieved
through various steps: we bound
$\Vert P_{j+2}\Vert_{{\cal L}(\cW_j,\cW_{j+2})}$ 
in section \ref{S:matrix-form} by 
the $\ell^\infty$-norm of a suitable matrix; in sections
\ref{S:first-bound} and \ref{S:second-bound} we characterize such an
$\ell^\infty$-norm and show that the desired bound reduces to certain
properties of a suitable function; we finally analyze such function
in section \ref{S:matrix-norm}. This analysis is computer assisted.}

\subsection{Bounding the operator norm by a matrix norm}\label{S:matrix-form}

Recalling the definition \eqref{eq:defWj} of the subspaces $\cW_j$ as well as Remark \ref{rem:parity}, we can split each $\cW_j$ into its two nontrivial orthogonal components according to parity; precisely, if $j$ is even we have $\cW_j= \cW_j^{(0,0)} \oplus \cW_j^{(1,1)}$, whereas if $j$ is odd we have $\cW_j= \cW_j^{(1,0)} \oplus \cW_j^{(0,1)}$. Furthermore, again by orthogonality it holds $P_{j+2} \in {\cal L}(\cW_j^{\, \alpha}, \cW_{j+2}^{\, \alpha})$ for any $\alpha \in \{0,1\}^2$; hence, our target result can be achieved by considering each parity component separately. Hereafter, we will analyze the case $j$ even and $\alpha =(0,0)$, i.e., we will bound $\Vert P_{j+2} \Vert_{{\cal L}(\cW_j^{(0,0)}, \cW_{j+2}^{(0,0)})}$; the other three cases can be treated similarly.

\smallskip
Let us set $d_j := \text{dim}\, \cW_j^{(0,0)}$, \rhn{note that $d_j = \frac{j}2-1$
because $j=|h|=h_1+h_2$ with $h_1,h_2\ge 2$ even}, and let us
introduce the normalized basis functions $\hat \Phi_h := \Phi_h /
\Vert \Phi_h \Vert_{H^1_0(\Omega)}$. For the sake of definiteness, let
us order the basis functions in each $\cW_j^{(0,0)}$ by increasing the
first index $h_1$. Any $v \in \cW_j^{(0,0)}$ is represented as
$$
v=\sideset{}{'} \sum_{|h|=j} \hat v_h \hat \Phi_h
\quad\textrm{with} \quad
{\mathbf v}:=(\hat v_h) \in \mathbb{R}^{d_j};
$$
the summation symbol means that only indices $h=(h_1,h_2)$ with even
components are considered. Similarly, any $w \in \cW_{j+2}^{(0,0)}$ is
represented as
$$
w=\sideset{}{'} \sum_{|k|=j+2} \hat w_k \hat \Phi_k
\quad\textrm{with} \quad
{\mathbf w}:=(\hat w_k) \in \mathbb{R}^{d_{j+2}}.
$$
Therefore, if $w=P_{j+2} v$, then ${\mathbf w}= {\mathbf A}^T_j {\mathbf v}$,
where ${\mathbf A}_j  \in \rhn{\mathbb{R}^{d_j\times d_{j+2}}}$ is the
matrix whose entries are
$$
a_{hk}:=(\hat \Phi_h,\hat \Phi_k)_{H^1_0(\Omega)}
\quad\textrm{for} \quad |h|=j, \ |k|=j+2,
$$
and  $h_1,\, h_2, \, k_1,\, k_2$ even (recall that the \rhn{$\Phi_h$'s} that span $\cW_j^{(0,0)}$ form an orthogonal basis for this space).
Note that ${\mathbf A}_j$ is a sub-block of the (even,\,even) block
\rhn{${\mathbf A}^{0,0}$} of the stiffness matrix ${\mathbf A}$ for
the normalized Babu\v ska-Shen basis in $H^1_0(\Omega)$. \rhn{More}
precisely, denoting by ${\mathbf I}_j \in \rhn{\mathbb{R}^{d_j \times
    d_j}}$ the identity matrix of order $d_j$, we have
\rhn{${\mathbf A}^{0,0} = {\sf tridiag} \ ({\mathbf A}^T_{j-2} , \
  {\mathbf I}_j , \ {\mathbf A}_j)$} (with $j \geq 4$).

Now, one immediately has 
$$
\Vert P_{j+2} \Vert_{{\cal L}(\cW_j^{(0,0)},\cW_{j+2}^{(0,0)})} = \Vert {\mathbf A}^T_j \Vert_2 = \Vert {\mathbf A}_j \Vert_2 = \Vert P_{j+2}^* \Vert_{{\cal L}(\cW_{j+2}^{(0,0)},\cW_j^{(0,0)})} 
$$
(where $\Vert \cdot \Vert_p$ denotes the $p$-norm of a matrix), which together with the inequality $\Vert {\mathbf A}_j \Vert_2 = 
\rho({\mathbf A}_j {\mathbf A}^T_j)^{1/2} \leq \Vert {\mathbf A}_j {\mathbf A}^T_j \Vert_\infty^{1/2}$, yields the bound 
\begin{equation}\label{eq:bound-projnorm}
\Vert P_{j+2} \Vert_{{\cal L}(\cW_j^{(0,0)},\cW_{j+2}^{(0,0)})} \leq \Vert {\mathbf A}_j {\mathbf A}^T_j \Vert_\infty^{1/2}\;.
\end{equation}

\subsection{A first expression for the matrix entries}\label{S:first-bound}

In order to compute the norm on the right-hand side \rhn{of
  \eqref{eq:bound-projnorm},} we observe that ${\mathbf A}_j$ is a
bi-diagonal matrix by condition \eqref{eq:A6}. \rhn{In fact, for any
index $h$ with $|h|=j$ the only indexes $h',h''$ with $|h'|=|h''|=j+2$ 
that give rise to entries $a_{h,h'}$ and $a_{h,h''}$ different from 0 are
$h'=h+(0,2)$ and $h''=h+(2,0)$.}
The explicit value of these entries is computable via the following formulas (in which all inner-products and norms are those of $H^1_0(\Omega)$):
$$
a_{h,h'}= (\hat \Phi_{h},\hat \Phi_{h'}) = \frac{(\Phi_h,\Phi_{h'})}{\Vert \Phi_h \Vert \,  \Vert \Phi_{h'} \Vert}
\qquad \qquad
a_{h,h''}= (\hat \Phi_h,\hat \Phi_{h''}) = \frac{(\Phi_h,\Phi_{h''})}{\Vert \Phi_h \Vert \, \Vert \Phi_{h''} \Vert}
$$ 
with
\begin{equation*}
\begin{split}
(\Phi_h,\Phi_{h'}) &=  -\frac1{(2h_2+1)\sqrt{(2h_2-1)(2h_2+3)}}, \\
(\Phi_h,\Phi_{h''}) &=-\frac1{(2h_1+1)\sqrt{(2h_1-1)(2h_1+3)}},  
\end{split}
\end{equation*}
and
$$
\Vert \Phi_h\Vert^2=\frac2{(2h_1-3)(2h_1+1)}+\frac2{(2h_2-3)(2h_2+1)},
$$
whence
$$
\Vert \Phi_{h'}\Vert^2=\frac2{(2h_1-3)(2h_1+1)}+\frac2{(2h_2+1)(2h_2+5)},
$$
$$
\Vert \Phi_{h''}\Vert^2=\frac2{(2h_1+1)(2h_1+5)}+\frac2{(2h_2-3)(2h_2+1)}.
$$

Since $2\leq h_1 \leq j-2$ ($h_1$ even), it is convenient to set $n :=\frac{j}2$ and $h_1 :=2i$, with $1 \leq i \leq n-1$; consequently,
$h_2:=j-h_1=2(n-i)$. Substituting these expressions in the previous formulas, we obtain
$$
(\Phi_h,\Phi_{h'})=  -\frac1{(4(n-i)+1)\sqrt{(4(n-i)-1)(4(n-i)+3)}} =: a_i
$$
$$
(\Phi_h,\Phi_{h''})=-\frac1{(4i+1)\sqrt{(4i-1)(4i+3)}} =: b_i,
$$
and
$$
\Vert \Phi_h\Vert^2=\frac2{(4i-3)(4i+1)}+\frac2{(4(n-i)-3)(4(n-i)+1)} =: \phi_i,
$$
$$
\Vert \Phi_{h'}\Vert^2=\frac2{(4i-3)(4i+1)}+\frac2{(4(n-i)+1)(4(n-i)+5)} =: \psi_i,
$$
$$
\Vert \Phi_{h''}\Vert^2=\frac2{(4i+1)(4i+5)}+\frac2{(4(n-i)-3)(4(n-i)+1)} =: \eta_i.
$$
Note that $b_i=a_{n-i}$ and $\phi_i=\phi_{n-i}$, $\eta_i=\psi_{n-i}$. Hence, for $1 \leq i \leq n-1$,
$$
({\mathbf A}_j)_{ii}= \frac{a_i}{\sqrt{\phi_i \psi_i}} =: \rhn{\delta_i,} \qquad 
({\mathbf A}_j)_{i,i+1}= \frac{b_i}{\sqrt{\phi_i \eta_i}} = \frac{a_{n-i}}{\sqrt{\phi_{n-i} \psi_{n-i}}} = \rhn{\delta_{n-i},}
$$
i.e., ${\mathbf A}_j = {\sf bidiag} \ \rhn{( \delta_i , \ \delta_{n-i}
  )}$. Consequently, the nonzero entries of the matrix ${\mathbf
  A}_j{\mathbf A}^T_j \in \rhn{{\mathbb R}^{d_j\times d_j}}$ are
$$
({\mathbf A}_j{\mathbf A}^T_j)_{i,i-1}= \rhn{\delta_i\delta_{n-i+1},}
\quad
({\mathbf A}_j{\mathbf A}^T_j)_{ii}= \rhn{\delta_i^2 + \delta_{n-i}^2,}
\quad
({\mathbf A}_j{\mathbf A}^T_j)_{i,i+1}= \rhn{\delta_{i+1}\delta_{n-i},}
$$
i.e., ${\mathbf A}_j{\mathbf A}^T_j= {\sf tridiag}
\ (\rhn{\delta_i\delta_{n-i+1},\ \delta_i^2 + \delta_{n-i}^2 , \ \delta_{i+1}\delta_{n-i} )}$. Let us denote by 
${\sf s}_i^{(j)}$ the sum
of the entries in the $i$-th row of the matrix ${\mathbf A}_j{\mathbf A}^T_j$, which are all non-negative. Setting for convenience
\rhn{$\delta_n=0$,} we thus have
\begin{equation}\label{eq:row-sum}
  {\sf s}_i^{(j)}=
\rhn{\delta_{i}\delta_{n-i+1} + \delta_i^2 + \delta_{n-i}^2 + \delta_{i+1}\delta_{n-i},}
\qquad 1 \leq i \leq n-1.
\end{equation}
It is easily seen that ${\sf s}_i^{(j)}={\sf s}_{n-i}^{(j)}$ for $1 \leq i \leq \frac{n}2$. Since
\begin{equation}\label{eq:bound-infnorm}
\Vert {\mathbf A}_j {\mathbf A}^T_j \Vert_\infty = \max_{1 \leq i \leq n-1} {\sf s}_i^{(j)} \;,
\end{equation}
in view of (\ref{eq:bound-projnorm}) we are left with the problem of proving the existence of a constant $C>0$ such that
\begin{equation}\label{eq:desired-bound}
\max_{1 \leq i \leq n-1} {\sf s}_i^{(j)} \leq \frac14  - \frac{C}{j^2} \qquad \text{ for all  \ } j \geq 4;
\end{equation}
indeed, thanks to $\sqrt{\frac{1}{4}-x}\leq \frac{1}{2} -x$ for $x\leq \frac{1}{4}$, we obtain Proposition~\ref{prop:proj} with $c=C$.

A direct computation shows that \rhn{${\sf s}_1^{(j)}$ and  ${\sf s}_{n-1}^{(j)}$}
satisfy the bound in \eqref{eq:desired-bound} for a suitable $C$,
\rhn{because ${\sf s}_1^{(j)}={\sf s}_{n-1}^{(j)} < \frac14$} for all $j \geq 4$ and ${\sf s}_1^{(j)} \to \frac3{28}$ as $j \to \infty$. Thus, in the sequel we focus on the rows indexed from 2 to $n-2$, for $j \geq 8$ (i.e., $n \geq 4$).

\subsection{A second expression for the matrix entries}\label{S:second-bound}

We now apply a change of variables. Observing that all quantities
$a_i$, $\phi_i$, $\psi_i$, \rhn{$\eta_i$, $\delta_i$} defined above
depend upon $4i$ or $4(n-i)$ \rhn{for $2\le i\le n-2$,} we first set
$I :=4i$ and $N :=4n \geq 16$.
\rhn{To introduce the new variables $(t,r)$, we first go back to the original
range $1\le i \le n-1$, i.e. $4 \leq I \leq N-4$, and parametrized
$I$ as follows
$$
I=4(1-t)+(N-4)t=4+Rt,
\qquad
0 \le t \le 1,
$$
with 
$R:=N-8 \geq 8$. Similarly, we write
$$
N-I = 4+R\tau,
\qquad
\tau=\tau(t) := 1-t.
$$
At last, we introduce the second} parameter $r:=\frac1R \leq \frac18$. With these notation at hand, we easily obtain the following expressions for  $a_i$, $\phi_i$ and $\psi_i$:
$$
a_i^2 = \frac1{R^4} \frac1{(\tau +3r)(\tau+5r)^2(\tau+7r)} =: \frac1{R^4} A(t, r),
$$
$$
\phi_i=\frac1{R^2} \left( \frac2{(t+r)(t+5r)}+\frac2{(\tau+r)(\tau+5r)} \right) =: \frac1{R^2} B(t,r),
$$
$$
\psi_i=\frac1{R^2} \left( \frac2{(t+r)(t+5r)}+\frac2{(\tau+5r)(\tau+9r)} \right) =: \frac1{R^2} C(t,r).
$$
Hence, we arrive at
$$
\rhn{\delta_i^2} = \frac{a_i^2}{\phi_i \psi_i} = \frac{A(t,r)}{B(t,r)C(t,r)} =: D(t,r).
$$
Straightforward computations show that
$$
\rhn{\delta_{i+1}^2} =D(t+4r,r),
\qquad
\rhn{\delta_{n-i}^2} =D(\tau,r),
\qquad
\rhn{\delta_{n-i+1}^2} =D(\tau+4r,r). 
$$
We conclude that the sum of the entries in the $i$-th row of ${\mathbf A}_j{\mathbf A}^T_j$, given by (\ref{eq:row-sum}), can be expressed as follows:
\rhn{
\begin{equation}\label{eq:def-Str}
\begin{aligned}
{\sf s}_i^{(j)} &= \sqrt{D(t,r)D(\tau+4r,r)} + D(t,r)
\\
& + \sqrt{D(t+4r,r)D(\tau,r)} + D(\tau,r) =: S(t,r)
\end{aligned}
\end{equation}
}
for  $2 \leq i \leq n-2$, which is equivalent to  $4r \leq t \leq 1-4r$.

\subsection{Bounding the matrix norm}\label{S:matrix-norm}

\begin{figure}[h!]
\begin{center}
\includegraphics[width=6.5cm]{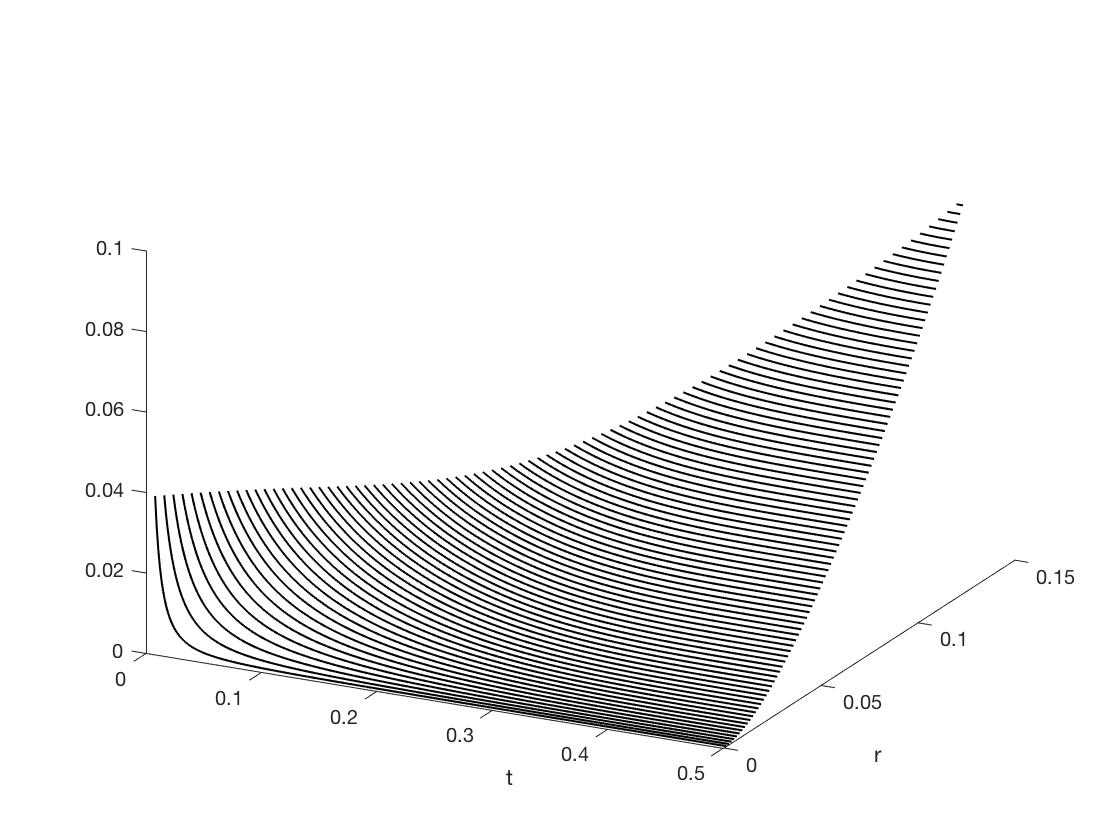}
\includegraphics[width=6.cm]{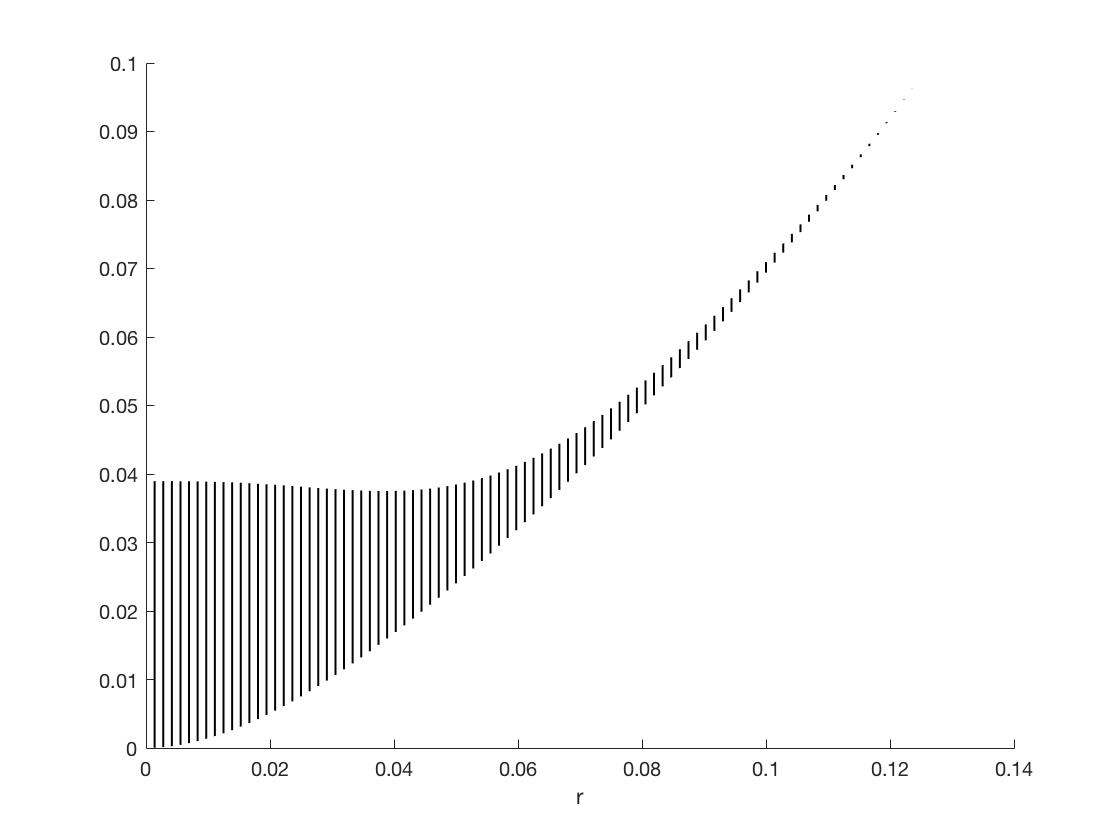}
\caption{Two views of the graph of the function \ $\frac14 - S(t,r)$
\rhn{on the triangle $0\le t\le \frac{1}{2}$, $0\le r \le
  \frac{1}{4}t$. Note that $S(t,r)$ is multi-valued at $t=r=0$.}}
\label{fig:plot_S}
\end{center}
\end{figure}

Since the function $S(t,r)$ is symmetric \rhn{with respect to} $t=\frac{1}{2}$ for any $r$, we may restrict it to the triangle $0 \leq t \leq \frac12$, 
$0 \leq r \leq \frac14 t$. Fig. \ref{fig:plot_S} \rhn{displays} two
plots of the function $\frac14 - S(t,r)$, \rhn{and suggests
clearly} that $S(t,r) < \frac14$ whenever $r>0$, with a quadratic behavior in $r$ at the origin.
However, \rhn{establishing such results rigorously is somehow complicated}
by the fact that $S(t,r)$ is singular at $(t,r)=(0,0)$, where it becomes multi-valued.

\begin{figure}[h!]
\begin{center}
\includegraphics[width=8cm]{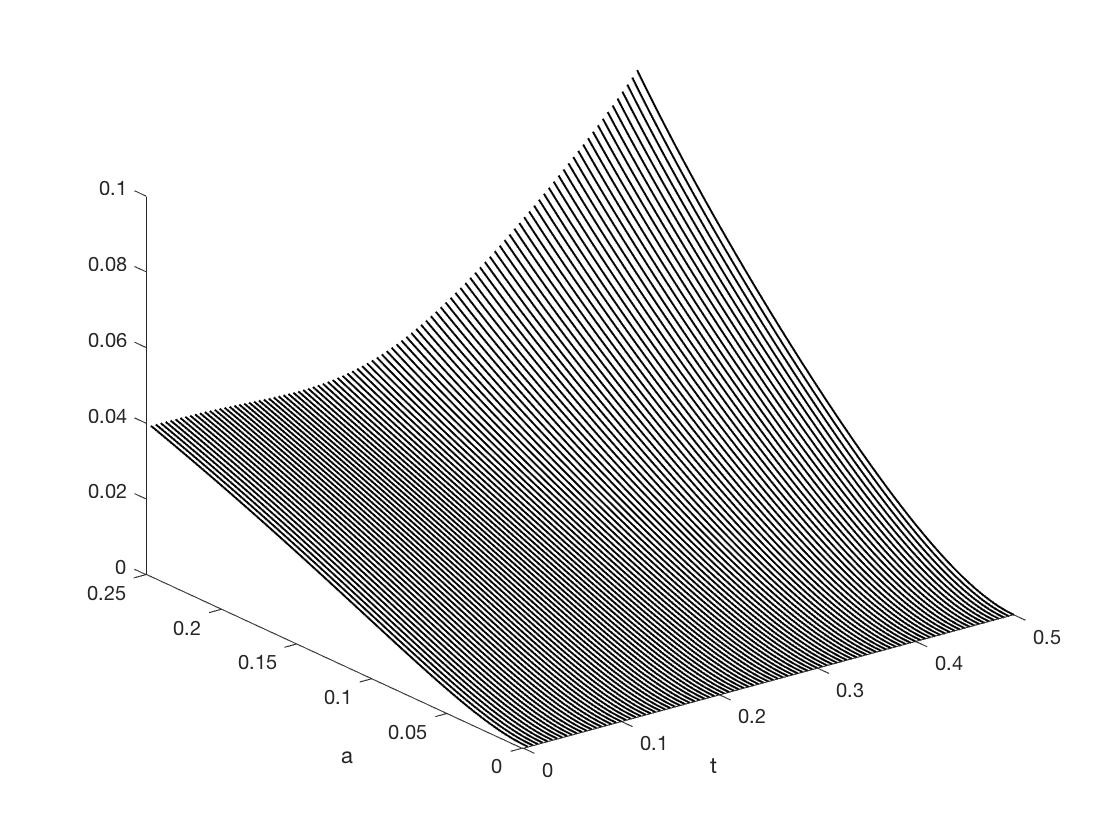}
\caption{A view of the graph of the function \ $\frac14 - \sigma(t,a)$
\rhn{in the rectangle $0\le t \le \frac{1}{2}, 0\le a \le \frac{1}{4}.$}}
\label{fig:plot_sigma}
\end{center}
\end{figure}

To remove this singularity, we apply the Duffy transform $(t,a) \mapsto (t,r)=(t,at)$, which maps the rectangle $0 \leq t \leq \frac12$, $0 \leq a \leq \frac14$ onto the triangle $0 \leq t \leq \frac12$, $0 \leq r \leq \frac14 t$.
Correspondingly, we are led to consider the function $\sigma(t,a) := S(t,at)$, which turns out to be smooth everywhere in this rectangle;
a plot of the function $\frac14 - \sigma(t,a)$ is \rhn{depicted in} Fig. \ref{fig:plot_sigma}.

With the help of a symbolic manipulator, we obtain, for all $t \in [0,\tfrac12]$,
$$
\sigma(t,0)=\frac14, \qquad \frac{\partial \sigma}{\partial a}(t,0) = 0, \qquad \frac{\partial^2 \sigma}{\partial a^2}(t,0) = -\frac{G(t)}{\tau^2},
$$
where
$$
G(t) := 3t^{10} + 9t^8\tau^2 - 8t^7\tau^3 + 16t^6\tau^4 + 24t^5\tau^5 + 16t^4\tau^6 - 8t^3\tau^7 + 9t^2\tau^8 + 3\tau^{10}.
$$
We note that the polynomial $G(t)$ is strictly decreasing in
$[0,\frac12]$ between $G(0)=3$ and $G(\frac12)=1$. \rhn{We thus easily see} that
$\frac{\partial^2 \sigma}{\partial a^2}(t,0) \leq -1$ for $0 \leq t \leq \frac12$; hence, by continuity we get the existence of two constants $C_*>0$ and  $a_* \in (0,\tfrac14]$ such that $\frac{\partial^2 \sigma}{\partial a^2}(t,a)  \leq  -C_*$ for $0 \leq t \leq \frac12$ and $0 \leq a \leq a_*$. With these constants at hand, by Taylor's expansion with Lagrange's reminder, we are entitled to write
$$
\sigma(t,a) \, = \, \frac14 + \frac{\partial^2 \sigma}{\partial a^2}(t,\bar{a}) a^2 \, \leq \, \frac14 - C_* a^2 
\qquad \text{for \ } 0\leq t \leq \tfrac12, \ 0 < a \leq a_*,
$$
with some $\bar{a} = \bar{a}(t,a) \in (0, a)$. 

\begin{figure}[h!]
\begin{center}
\includegraphics[width=11cm]{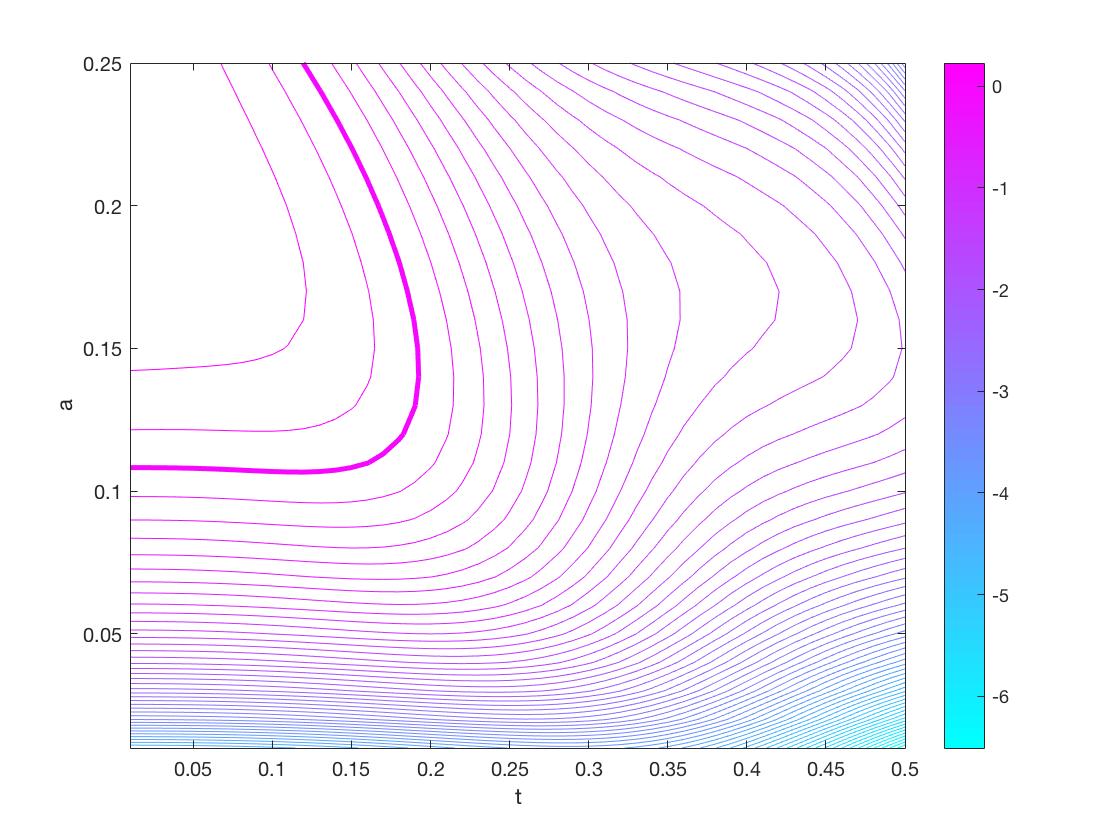}
\caption{Contour plot of  $\frac{\partial^2 \sigma}{\partial a^2}(t,a)$: the function is negative on $[0,\frac12]\times [0,a^*]$ with 
$a^*>\frac{1}{10}$. The thicker curve represents the zero level set}
\label{fig:plot_der_2_sigma}
\end{center}
\end{figure}

By computing the symbolic expression of the function $\frac{\partial^2
  \sigma}{\partial a^2}(t,a)$ and by \rhn{examining} its level sets
(via a numerical procedure), one finds that $a_*$ safely satisfies $a_* > \frac1{10}$ (see Figure \ref{fig:plot_der_2_sigma}). 
Therefore, going back to our function $S(t,r)=\sigma(t,\frac{r}{t})$, we deduce that
$$
S(t,r) \leq \frac14 - \frac{C_*}{t^2} r^2 \leq \frac14 - 4C_* r^2 \qquad  \text{for \ } 0 < r \leq \tfrac1{10}t, \ \  t \leq \tfrac12.
$$
Recalling \eqref{eq:def-Str} and using the expressions $t=4r(i-1)$ and $r=\frac12\, \frac1{j-4} > \frac1{2j}$, we immediately obtain
\begin{equation}\label{eq:bound-sij}
{\sf s}_i^{(j)} \leq \frac14 - \frac{C_*}{j^2} \qquad  \text{for \ } 4 \leq i \leq \tfrac{n}2.
\end{equation}
\rhn{Note that we require the restriction $i\ge 4$ to satisfy the constraint
$$
  r \le \frac{1}{10} t = \frac{2}{5} r(i-1)
  \qquad\Rightarrow\qquad
  i \ge \frac{7}{2}.
$$
Therefore, we are left with the task of establishing a similar bound for
${\sf s}_2^{(j)}$ and ${\sf s}_3^{(j)}$ by different means.}
It is easily checked that  for $j \to \infty$ it holds ${\sf s}_2^{(j)} \to \frac{65}{308} < \frac14$ and ${\sf s}_3^{(j)} \to \frac{55}{220} < \frac14$, while both ${\sf s}_2^{(j)}$ and ${\sf s}_2^{(j)}$ are $<\frac14$ for all $j \geq 8$. This implies the desired bound for a suitable constant $C_{**}>0$.   

\smallskip
The proof of \eqref{eq:desired-bound} is thus \rhn{complete,
whence} Proposition \ref{prop:proj} is established.
\endproof

\section*{Acknowledgements}
The authors wish to thank Valeria Simoncini for helpful suggestions concerning Section 6. The first and fourth authors are members of the INdAM research group GNCS, which granted partial support to this research. The second author has been artially supported by the NSF Grant DMS -1411808, the Institut Henri Poincar\'e (Paris) and the Hausdorff Institute (Bonn).



\begin{thebibliography}{99}
  
\bibitem{Babuska-1991}
I.~Babu{\v{s}}ka, A.~Craig, J.~Mandel, and J.~Pitk{\"a}ranta.
\newblock Efficient preconditioning for the p-version finite element method in
  two dimensions.
  \newblock {\em {SIAM J. Numer. Anal.}}, 28(3):624Ð--661, 1991.

\bibitem{BW85}
R.E. Bank and A. Weiser.
\newblock  Some a posteriori error estimators for elliptic partial 
differential equations
\newblock {\em Math. Comp.}, 44:285--301, 1985.

\bibitem{BEK:96}
F. A. Bornemann, B. Erdmann, and R. Kornhuber,
\newblock A posteriori error estimates for elliptic problems in
	two and three space dimensions,
\newblock  \emph{SIAM J. Numer. Anal.}, 33:1188--1204, 1996.
        
\bibitem{69.1}
M.~B{\"u}rg and W.~D{\"o}rfler.
\newblock Convergence of an adaptive {$hp$} finite element strategy in higher
  space-dimensions.
\newblock {\em Appl. Numer. Math.}, 61(11):1132--1146, 2011.

\bibitem{21}
P.~Binev, W.~Dahmen, and R.~{DeV}ore.
\newblock Adaptive finite element methods with convergence rates.
\newblock {\em Numer. Math.}, 97(2):219 -- 268, 2004.

\bibitem{Binev-2015}
P.~Binev.
\newblock Tree approximation for $hp$-adaptivity.
\newblock IMI Preprint Series 2015:07, 2015.

\bibitem{BPS-2009}
D.~Braess, V.~Pillwein, and J.~Sch{\"o}berl.
\newblock Equilibrated residual error estimates are {$p$}-robust.
\newblock {\em Comput. Methods Appl. Mech. Engrg.}, 198(13-14):1189--1197,
  2009.

\bibitem{CHQZ06}
C.~Canuto, M.~Y.Hussaini, A.~Quarteroni, and T.~A.~Zang.
\newblock {\sl Spectral Methods. Fundamentals in Single Domains.}
\newblock Springer 2016.

\bibitem{CNSV16} 
C.~Canuto, R.H.~Nochetto, R. Stevenson, and M.~Verani.
\newblock Convergence and optimality of $hp$-AFEM.
\newblock \emph{Numer. Math.} 135 (2017), 1073-1119

\bibitem{CNSV17} 
C.~Canuto, R.H.~Nochetto, R. Stevenson, and M.~Verani.
\newblock On $p$-robust saturation for $hp$-AFEM.
\newblock \emph{Comput. \& Math. with Appl. } 73 (2017),  2004--2022.

\bibitem{64.145}
L.~Demkowicz, J.~Gopalakrishnan, and J.~Sch{\"o}berl.
\newblock Polynomial extension operators. {P}art {III}.
\newblock {\em Math. Comp.}, 81(279):1289--1326, 2012.

\bibitem{DN02}
W. D\"orfler and R.H. Nochetto, 
\newblock Small data oscillation implies the saturation assumption.
\newblock \emph{Numer. Math.} 91:1--12, 2002.

\bibitem{ErnVohralik-2015}
A.~Ern and M.~Vohral{\'{\i}}k.
\newblock Polynomial-degree-robust a posteriori estimates in a unified setting
  for conforming, nonconforming, discontinuous {G}alerkin, and mixed
  discretizations.
\newblock {\em SIAM J. Numer. Anal.}, 53(2):1058--1081, 2015.

\bibitem{ErnVohralik-2016}
A.~Ern and M.~Vohral{\'{\i}}k.
\newblock Stable broken {H1} and {H(div)} polynomial extensions for
  polynomial-degree-robust potential and flux reconstruction in three space
  dimensions.
\newblock hal 01422204, {Inria Paris-Rocquencourt}, 2016.

\bibitem{GuoBabuska-1}
B.~Guo and I.~Babu{\v{s}}ka.
\newblock {The $hp$ version of finite element method, Part 1: The basic
  approximation results}.
\newblock {\em {Comp. Mech}}, 1:21--41, 1986.

\bibitem{GuoBabuska-2}
B.~Guo and I.~Babu{\v{s}}ka.
\newblock {The $hp$ version of finite element method, Part 2: General results
  and application}.
\newblock {\em {Comp. Mech}}, 1:203--220, 1986.

\bibitem{MW01}
J.~M. Melenk and B.~I. Wohlmuth.
\newblock On residual-based a posteriori error estimation in {$hp$}-{FEM}.
\newblock {\it Adv. Comput. Math.}, 15(1-4):311--331, 2002.

\bibitem{Noch93}
R.H. Nochetto,
\newblock Removing the saturation assumption in a posteriori error
analysis,
\newblock  \emph{Istit. Lombardo Accad. Sci. Lett. Rend. A},
127:67-82, 1993.

\end{thebibliography}
\end{document}